# Strong Low Degree Hardness for the Number Partitioning Problem

Rushil Mallarapu[*]  Mark Sellke[†]

**Abstract**

In the *number partitioning problem (NPP)* one aims to partition a given set of $N$ real numbers into two subsets with approximately equal sum. The NPP is a well-studied optimization problem and is famous for possessing a *statistical-to-computational gap*: when the $N$ numbers to be partitioned are i.i.d. standard gaussian, the optimal discrepancy is $2^{-\Theta(N)}$ with high probability, but the best known polynomial-time algorithms only find solutions with a discrepancy of $2^{-\Theta(\log^2 N)}$. This gap is a common feature in optimization problems over random combinatorial structures, and indicates the need for a study that goes beyond worst-case analysis.

We provide evidence of a nearly tight algorithmic barrier for the number partitioning problem. Namely we consider the family of low coordinate degree algorithms (with randomized rounding into the Boolean cube), and show that degree $D$ algorithms fail to solve the NPP to accuracy beyond $2^{-\tilde{O}(D)}$. According to the low degree heuristic, this suggests that simple brute-force search algorithms are nearly unimprovable, given any allotted runtime between polynomial and exponential in $N$. Our proof combines the isolation of solutions in the landscape with a conditional form of the overlap gap property: given a good solution to an NPP instance, slightly noising the NPP instance typically leaves no good solutions near the original one. In fact our analysis applies whenever the $N$ numbers to be partitioned are independent with uniformly bounded density.

## Contents

---
[*]Harvard University. Email: rushil.mallarapu@pm.me.
[†]Department of Statistics, Harvard University. Email: msellke@fas.harvard.edu.

# 1 Introduction

The number partitioning problem (NPP) asks to partition $N$ given real numbers into two parts with approximately equal sum. That is, given real numbers $g_1, ..., g_N$, one aims to find the subset $A$ of $[N] := \{1, 2, ..., N\}$ which minimizes the discrepancy

$$\left| \sum_{i \in A} g_i - \sum_{i \notin A} g_i \right|.$$

Identifying the input $g_1, ..., g_N$ with a vector $g \in \mathbf{R}^N$, this is equivalent to choosing a point $x$ in the $N$-dimensional binary hypercube $\Sigma_N := \{\pm 1\}^N$ which minimizes the discrepancy:

$$\min_{x \in \Sigma_N} |\langle g, x \rangle|. \tag{1.1}$$

Rephrased as a decision problem – whether there exists a subset $A \subseteq [N]$ (or a point $x \in \Sigma_N$) such that the discrepancy is zero, or sufficiently small – the NPP is NP-complete as can be shown by reduction from subset sum. In fact, the NPP is one of the six basic NP-complete problems of Garey and Johnson, and of those, the only one involving numbers [GJ79, § 3.1].

Solving the NPP has a number of practical applications. Together with a multiway generalization, it was first formulated by Graham, who considered it in the context of multiprocessor scheduling [Gra69]. Later work by Coffman, Garey, and Johnson, as well as Tsai, looked at utilizing algorithms for the NPP for designing multiprocessor schedulers or large integrated circuits [CGJ78], [Tsa92]. Coffman and Lueker further discuss how the NPP can be used for resource allocation problems [CL91]. The NPP was also used to design an early public key cryptography scheme [MH78], later shown to be insecure by [Sha82].

An important application of the NPP to statistics comes from the design of *randomized controlled trials* [KAK19], [Har+24]. Consider $N$ individuals, each with a covariate vector $g_i \in \mathbf{R}^d$. Then the problem is to divide them into a treatment group (denoted $A_+$) and a control group (denoted $A_-$), while ensuring that the covariates across both groups are balanced. In our notation, this amounts to finding an $A_+$ (with $A_- := [N] \setminus A_+$) minimizing

$$\min_{A_+ \subseteq [N]} \left\| \sum_{i \in A_+} g_i - \sum_{i \in A_-} g_i \right\|_\infty. \tag{1.2}$$

This vector extension of the NPP (the $d = 1$ case) is called the vector balancing problem (VBP).

Another major source of interest in the NPP, as well as potential explanations for its algorithmic hardness, comes from statistical physics. In the 1980s, Derrida introduced the eponymous *random energy model (REM)*, a simplified example of a spin glass in which, unlike e.g. the Sherrington-Kirkpatrick model, the energy levels of different states are independent of each other [Der80], [Der81], [BFM04]. Despite this simplicity, this model enabled heuristic analyses of the Parisi theory for mean field spin glasses, and it was suspected that very general disordered systems would locally behave like the REM and its generalizations [BM04], [Kis15]. The NPP was the first system for which this local REM conjecture was confirmed [Bor+09a], [Bor+09b]. Relatedly, in the case when the $g_i$ are chosen independently and uniformly from the discrete set $\{1, 2, 3, ..., 2^M\}$, Gent and Walsh conjectured that the hardness of finding perfect partitions (i.e., with discrepancy zero if $\sum_i g_i$ is even, and one else) was controlled by the parameter $\kappa := \frac{M}{N}$ [GW98], [GW00]. Mertens soon gave a nonrigorous statistical mechanics argument suggesting the existence of a phase transition from $\kappa < 1$ to $\kappa > 1$; that is, while solutions exist in the low $\kappa$ regime, they stop existing in the high $\kappa$ regime [Mer01]. It was also observed that this phase transition coincided with the empirical onset of computational hardness for typical algorithms, and Borgs, Chayes, and Pittel proved the existence of this phase transition soon after [Hay02], [BCP01].

**The Statistical-to-Computational Gap.** Many problems involving searches over random combinatorial structures (i.e., throughout high-dimensional statistics) exhibit a statistical-to-computational gap: the optimal values which are known to exist via non-constructive, probabilistic methods are far better than those achievable by state-of-the-art algorithms. In the pure optimization setting, examples of such gaps are found in random constraint satisfaction [MMZ05], [AC08], [Kot+17], finding maximal independent sets in sparse random graphs [GS14], [CE15], the largest submatrix problem [GL18], [GJS21], and the $p$-spin and diluted $p$-spin models [GJ21], [Che+19]. These gaps also arise in various planted models, such as matrix or tensor PCA [BR13], [LKZ15a], [LKZ15b], [HSS15], [Hop+17], [AGJ20], high-dimensional regression [GZ22], or the infamously hard planted clique problem [Jer92], [DM15], [MPW15], [Bar+19], [GZ24]. These indicate that these problems are "hard" beyond being in NP: algorithms fail even for average-case instances.

The NPP is no exception, and exhibits an exponentially wide statistical-to-computational gap. For the global optimum, Karmarkar-Karp-Lueker-Odlyzko showed that when the $g_i$ are i.i.d. random variables with sufficiently nice distribution,[1] the minimum discrepancy of (1.1) is $\Theta\bigl(\sqrt{N} 2^{-N}\bigr)$ with high probability as $N \to \infty$ [Kar+86]. Their result also extends to *even partitions*, where the sizes of each subset is equal (i.e., for $N$ even), worsening only to $\Theta(N 2^{-N})$. Yet the best known algorithms cannot achieve discrepancies anywhere close to this in polynomial time. A simple greedy algorithm based on sorting $g_1, ..., g_N$ and grouping adjacent pairs is known to achieve discrepancy $O(N^{-1})$ (see [Mer06]), and more recently [KAK19] gave another algorithm achieving discrepancy $O(N^{-2})$ for balanced partitions. The best known algorithm is however due to Karmarkar-Karp [KK83] from 1983, and achieves discrepancy $O\bigl(N^{-O(\log N)}\bigr) = 2^{-O(\log^2 N)}$. See also [Lue87] and [Yak96] for performance bounds on the related *partial differencing method (PDM)* and *largest differencing method (LDM)*, and [BM08] for a conjecturally sharper description of the performance of the latter. The situation for the multidimensional VBP is similar: restricting for simplicity to fixed $d$ as $N$ grows, the minimum discrepancy is typically $\Theta\bigl(\sqrt{N} 2^{-N/d}\bigr)$ [TMR20], with an extension of the Karmarkar-Karp algorithm achieving discrepancy $2^{-\Theta(\log^2 N)/d}$.

**Algorithmic Hardness and Landscape Obstructions.** A broadly successful approach to random optimization problems has emerged recently based on analyzing the geometry of the solution landscape. Many "hard" random optimization problems have a certain disconnectivity property, in which solutions tend to live in isolated separated clusters [MMZ05], [ACR11], [AC08], [AMS25], [GJK23] . In the seminal work [GS14], Gamarnik and Sudan showed how to use a strong form of this disconnectivity to deduce rigorous hardness results for classes of stable algorithms, in what has become known as the *overlap gap property* (OGP) [Gam21]. In its simplest form, an OGP states that for certain constants $0 \leq \nu_1 < \nu_2$, for every two near-optimal states $x, x'$ for a particular instance $g$ of the problem either have $d(x, x') < \nu_1$ or $d(x, x') > \nu_2$. That is, pairs of solutions are either close to each other, or much further away – the condition that $\nu_1 < \nu_2$ ensures that the "diameter" of solution clusters is smaller than the separation between these clusters. The original work [GS14] studied maximal independent sets in sparse random graphs, obtaining suboptimal bounds later sharpened by [RV17], [Wei22]. However the method has been greatly developed since then, finding applications to a host of other problems such as largest average submatrix [GL18], mean-field spin glasses and related combinatorial optimization problems [GJ21], [GJW24], [Che+19], [HS25a], [HS23], [Jon+23], graph alignment [DGH25], (Not-All-Equal)-$k$-SAT [GS17], [BH21], and the random perceptron [Gam+22], [Gam+23], [LSZ24].

---

[1]Specifically, having bounded density and finite fourth moment.

For the NPP, it was expected for decades that the "brittleness" of the solution landscape would be a central barrier in finding successful algorithms to close the statistical-to-computational gap. Mertens wrote in 2001 that any local heuristics, which only looked at fractions of the domain, would fail to outperform random search [Mer01, § 4.3]. This was backed up by the failure of many algorithms based on locally refining Karmarkar-Karp-optimal solutions, such as simulated annealing [AFG96], [SFD96], [Joh+89], [Joh+91], [Ali+05]. Previously, Gamarnik and Kızıldağ applied the OGP methodology to the NPP, proving that a more complicated multi-OGP holds for discrepancies of $2^{-\Theta(N)}$ (i.e., the statistical near-optimum), but is absent for smaller discrepancies of $2^{-E_N}$ with $\omega(1) \leq E_N \leq o(N)$ [GK23]. However they were able to prove that for $\varepsilon \in (0, 1/5)$, no stable algorithm (suitably defined) could find solutions with discrepancy $2^{-E_N}$ for $\omega\left(N \log^{-\frac{1}{5}+\varepsilon} N\right) \leq E_N \leq o(N)$ [GK23, Thm. 3.2]. These results point to the efficacy of using landscape obstructions to show algorithmic hardness for the NPP, which we will take advantage of in Section 2.

Other evidence of hardness for the NPP has been provided by [Hob+17], who showed that a polynomial-time approximation oracle achieving discrepancy $O\left(2^{\sqrt{N}-N}\right)$ could give polynomial-time approximations in Minkowski's theorem on symmetric convex bodies. Very recently, Vafa and Vaikuntanathan showed that the Karmarkar-Karp algorithm's performance is nearly optimal among polynomial time algorithms, assuming the worst-case hardness of the shortest vector problem on lattices [VV25].

## 1.1 Statements of Main Results

We will use the following convenient terminology for the quality of an NPP solution.

**Definition 1.1.** Let $g \in \mathbf{R}^N$ be an instance of the NPP, and let $x \in \Sigma_N$. The *energy* of $x$ is

$$E(x; g) := -\log_2|\langle g, x \rangle|.$$

The *solution set* $S(E; g)$ is the set of all $x \in \Sigma_N$ that have energy at least $E$, i.e., that satisfy

$$|\langle g, x \rangle| \leq 2^{-E}.$$

Observe here that minimizing the discrepancy $|\langle g, x \rangle|$ corresponds to maximizing the energy $E(x; g)$; this terminology is motivated by the statistical physics literature [Mer01]. Rephrasing the preceding discussion, the *statistically optimal energy level* is $E = \Theta(N)$, while the best *computational energy level* currently known to be achievable in polynomial time is $E = \Theta(\log^2 N)$.

For our purposes, a *randomized algorithm* is a measurable function $\mathcal{A}: (g, \omega) \mapsto x \in \Sigma_N$ where $\omega \in \Omega_N$ is an independent random variable in some Polish space $\Omega_N$. Such an $\mathcal{A}$ is *deterministic* if it does not depend on $\omega$. For our main analysis, considering deterministic $\Sigma_N$-valued algorithms will suffice. We consider two classes of *low degree algorithms*, given by either low degree polynomials or by functions with low *coordinate degree*. Their study is motivated by the well-established *low degree heuristic*: degree $D$ algorithms (in either sense) are believed to be good proxies for the broader class of $e^{\tilde{O}(D)}$-time algorithms [Hop18], [KWB19]. At the same time, these classes are known to enjoy favorable stability properties making them amenable to rigorous analysis.

Our results show *strong low degree hardness* for the NPP at essentially all energy levels between the statistical and computational thresholds, in the sense of [HS25b].

**Definition 1.2** (Strong Low Degree Hardness [HS25b, Def. 3]). A sequence of random search problems, that is, a $N$-indexed sequence of random input vectors $g_N \in \mathbf{R}^{d_N}$ and random subsets

$$S_N = S_N(g_N) \subseteq \Sigma_N,$$

exhibits *strong low degree hardness (SLDH) up to degree* $D \le o(D_N)$ if, for all sequences of degree $o(D_N)$ algorithms $\mathcal{A}_N \colon (g, \omega) \mapsto x$ with $\mathbf{E}\|\mathcal{A}(y_N)\|^2 \le O(N)$, we have

$$\mathbf{P}(\mathcal{A}_N(g_N, \omega) \in S_N) \le o(1).$$

There are two related notions of degree which we want to consider in Definition 1.2. The first is traditional polynomial degree, applicable for algorithms given in each coordinate by low degree polynomial functions of the inputs. In this case, we show:

**Theorem 1.3** (Hardness for LDP Algorithms). *Let $g \sim \mathcal{N}(0, I_N)$ be a standard Normal random vector. Let $\mathcal{A}$ be any $\Sigma_N$-valued algorithm with $\mathbf{E}\|\mathcal{A}(g, \omega)\|^2 \le CN$, which obeys the stability estimate (1.5) below with degree $D$ in the LDP case. Assume $\omega(\log N) \le E \le N$ and $D \le o(2^{E/4})$. Then $\mathbf{P}(\mathcal{A}(g, \omega) \in S(E; g)) = o(1)$.*

Let us make a few comments about this result. First, the phrasing does not actually require $\mathcal{A}$ to be a polynomial, because a finite-degree polynomial cannot actually map standard Gaussian inputs to $\Sigma_N$. However this aside, the implications of Theorem 1.3 under the low degree heuristic are unreasonably pessimistic and primarily indicate that low degree polynomials are a poor proxy for efficient algorithms. For example, Theorem 1.3 rules out $N^{\Omega(\log N)}$ polynomial algorithms from even achieving the same energy threshold as the (polynomial time) Karmarkar-Karp algorithm, but the low degree heuristic suggests that such polynomials constitute a good proxy for all exponential time algorithms. On the other extreme, the low degree heuristic suggests that polynomial algorithms require doubly exponential time to achieve the statistically optimal discrepancy $E = \Theta(N)$, while brute-force search in fact requires at most exponential time.

Thus we turn to the second, more expressive notion of *coordinate degree*: a function $f \colon \mathbf{R}^N \to \mathbf{R}$ has coordinate degree $D$ if it can be expressed as a linear combination of functions depending on combinations of no more than $D$ coordinates. As coordinate degree is always smaller than polynomial degree, this enables us to consider a far broader class of *low coordinate degree (LCD)* algorithms. In fact we will consider a slightly broader class of $\mathbf{R}^N$-valued LCD algorithms, supplemented by an additional rounding scheme into $\Sigma_N$. The randomized rounding $\tilde{\mathcal{A}}$ of $\mathcal{A}$ is defined by $\tilde{\mathcal{A}} \coloneqq \mathsf{round}_{\tilde{U}}(\mathcal{A})$. Here $\mathsf{round}_{\tilde{U}}(x_1, ..., x_n) = \big(\mathsf{round}_{U_1}(x_1), ..., \mathsf{round}_{U_N}(x_N)\big) = (\mathsf{sign}(x_1 - U_1), ..., \mathsf{sign}(x_N - U_N))$ for i.i.d. Uniforms $U_i \in [-1, 1]$ and

$$\text{sign}(z) := \begin{cases} +1 & z > 0, \\ -1 & z \leq 0. \end{cases}$$

**Theorem 1.4** (Hardness for LCD Algorithms)**.** *Let $g \sim \mathcal{N}(0, I_N)$ be a standard Normal random vector. Let $\mathcal{A}$ be any coordinate degree $D$ algorithm with $\mathbf{E}\|\mathcal{A}(g)\|^2 \leq CN$, and let $\tilde{\mathcal{A}}$ be its randomized rounding as described above. Assume that*
   (a) *if $E = \delta N$ for $\delta \in (0,1)$, then $D \leq o(N)$;*
   (b) *if $\omega(\log^2 N) \leq E \leq o(N)$, then $D \leq o(E/\log^2(N/E))$.*

*Then $\mathbf{P}\big(\tilde{\mathcal{A}} \in S(E; g)\big) = o(1)$.*

**Remark 1.5.** Although Theorem 1.4 is stated for Gaussian disorder for simplicity and consistency with the low degree polynomial case, the proof extends verbatim to the much more general situation that $(g_1, ..., g_N)$ are independent random variables with uniformly bounded density (with respect to Lebesgue measure). Note that if each $g_i$ has density uniformly at most some constant $L$ independent of $N$, then the same holds (with the same value of $L$) for any independent sum $\sum_{s \in S} g_s$ for deterministic $S \subseteq [N]$. This is the only property of the Gaussian density that is needed to prove Theorem 1.4, appearing in Lemma 2.1 (which is used to show Lemma 2.4) and Lemma 2.12.

## 1.2 Heuristic Optimality of Theorem 1.4

Theorem 1.4 is essentially best-possible under the low degree heuristic. In particular, this heuristic suggests from the energy-degree tradeoff of $D \leq \tilde{o}(E)$ that finding solutions with energy $E$ requires time $e^{\tilde{\Omega}(E)}$. This tradeoff is attainable along the full range $1 \ll E \leq N$ via the following restricted brute-force search:
   (a) Choose a subset $J \subseteq [N]$ of $E$ coordinates (say, the first $E$).
   (b) Run an existing NPP algorithm on the restricted instance $g_{\bar{J}}$ to find $x_{\bar{J}}$ with $\langle g_{\bar{J}}, x_{\bar{J}} \rangle \leq 1$.
   (c) Fixing $x_{\bar{J}}$, the NPP given by $g$ turns into finding $x_J$ minimizing $|\langle g, x\rangle| = |\langle g_J, x_J\rangle + \langle g_{\bar{J}}, x_{\bar{J}}\rangle|$.
   (d) Fixing $j \in J$ and the coordinate $x_j = 1$, the conditional law of $\langle g_{\bar{J}}, x_{\bar{J}}\rangle + g_j$ is contiguous with a standard Gaussian. If $\langle g_{\bar{J}}, x_{\bar{J}}\rangle + g_j$ were exactly standard Gaussian, then by the previously mentioned results on statistically optimal solutions, with high probability there would exist $g$ attaining energy $\Theta(E)$. By contiguity, the same holds even with the shift $\langle g_{\bar{J}}, x_{\bar{J}}\rangle$.
   (e) Given the previous step, we can brute force search over the remaining $E$ coordinates in $J$, yielding a solution with energy $\Theta(E)$ with high probability in $e^{O(E)}$ time.

It is reasonable to ask whether the low (coordinate) degree heuristic is truly appropriate in our setting. In previous problems where it has been applied, the objective under consideration is stable to small input perturbations. By contrast solutions to the NPP are quite brittle, a fact which guides our proof. This brittleness also underlies the failure of low-degree polynomials in Theorem 1.3, which stems from the fact that polynomials can depend on the entries $g_i$ only in a somewhat coarse way. Nevertheless the sharpness of Theorem 1.4 is notable and highly suggestive. See also [LS24] for related discussion. We also mention that conjecturally optimal trade-offs of a similar flavor between running time and signal-to-noise ratio have been established for sparse PCA in [AWZ23], [Din+24] and for tensor PCA in [KWB19], [WEM19].

Finally, we note that algorithms with coordinate degree $\Omega(N)$ by definition involve nonlinear interactions between a constant fraction of the $N$ coordinates. Thus Theorem 1.4 essentially states that good NPP algorithms must be "truly global", which is echoed by recent heuristic algorithms for the NPP [KKS14], [COY19], [SBD21].

## 1.3 Notations and Preliminaries

We use the standard Bachmann-Landau notations $o(\cdot), O(\cdot), \omega(\cdot), \Omega(\cdot), \Theta(\cdot)$, in the limit $N \to \infty$. We abbreviate $f(N) \ll g(N)$ or $f(N) \gg g(N)$ when $f(N) = o(g(N))$ or $f(N) = \omega(g(N))$. In addition, we write $f(N) \lesssim g(N)$ or $f(N) \gtrsim g(N)$ when there exists an $N$-independent constant $C$ such that $f(N) \leq Cg(N)$ or $f(N) \geq Cg(N)$ for all $N$.

We write $[N] := \{1, ..., N\}$. If $S \subseteq [N]$, then we write $\overline{S} := [N] \setminus S$ for the complimentary set of indices. If $x \in \mathbf{R}^N$ and $S \subseteq [N]$, then the *restriction of $x$ to the coordinates in $S$* is the vector $x_S$ with

$$(x_S)_i := \begin{cases} x_i & i \in S, \\ 0 & \text{else.} \end{cases}$$

In particular, for $x, y \in \mathbf{R}^N$, $\langle x_S, y \rangle = \langle x, y_S \rangle = \langle x_S, y_S \rangle$.

On $\mathbf{R}^N$, we write $\|\cdot\|$ for the Euclidean norm, and $B(x, r) := \{y \in \mathbf{R}^N : \|y - x\| < r\}$ for the Euclidean ball of radius $r$ around $x$. In addition, we write

$$B_\Sigma(x, r) := B(x, r) \cap \Sigma_N = \{y \in \Sigma_N : \|y - x\| < r\},$$

to denote points on $\Sigma_N$ within distance $r$ of $x$.

We use $\mathcal{N}(\mu, \sigma^2)$ to denote the scalar Normal distribution with given mean and variance, and the abbreviation "r.v." for random variable (or random vector, if it is clear from context).

**$p$-Correlated and $p$-Resampled Pairs.** For $p \in [0, 1]$ and a pair $(g, g')$ of $N$-dimensional standard Normal random vectors, we say $(g, g')$ are *$p$-correlated* if $g'$ is distributed (conditionally on $g$) as

$$g' = pg + \sqrt{1 - p^2}\tilde{g},$$

where $\tilde{g}$ is an independent copy of $g$. We say $(g, g')$ are *$p$-resampled* if $g$ is a standard Normal random vector and $g'$ is drawn as follows: for each $i \in [N]$ independently,

$$g'_i = \begin{cases} g_i & \text{with probability } p, \\ \text{drawn from } \mathcal{N}(0, 1) & \text{with probability } 1 - p. \end{cases}$$

We denote such a pair by $g' \sim \mathcal{L}_p(g)$. In both cases, $g$ and $g'$ are marginally multivariate standard Normal and have entrywise correlation $p$. However, while

$$\mathbf{P}_{g' \sim_p g}(g = g') = 0,$$

when $(g, g')$ are $p$-resampled, we have

$$\mathbf{P}_{g' \sim \mathcal{L}_p(g)}(g = g') = \prod_{i=1}^N \mathbf{P}(g_i = g_{i'}) = (1 - \varepsilon)^N. \tag{1.3}$$

which is non-negligible for $\varepsilon \leq O(1/N)$.

**Coordinate Degree.** Let $\gamma_N$ be the $N$-dimensional standard Normal measure on $\mathbf{R}^N$, and consider the space $L^2(\gamma_N)$; this is the space of $L^2$ functions of $N$ i.i.d. standard Normal r.v.s. For $g \in \mathbf{R}^N$ and $S \subseteq [N]$, we can define subspaces of $L^2(\gamma_N)$

$$V_S := \{f \in L^2(\gamma_N) : f(g) \text{ depends only on } g_S\},$$
$$V_{\leq D} := \mathrm{span}(\{V_J : J \subseteq [N], |J| \leq D\})$$

These subsets describe functions which only depend on some subset of coordinates, or on some bounded number of coordinates. Note that $V_{[N]} = V_{\leq N} = L^2(\mathbf{R}^N, \pi^{\otimes N})$. The *coordinate degree* of a function $f \in L^2(\gamma_N)$ is defined as $\min\{D : f \in V_{\leq D}\}$. Note that if $f$ is a degree $D$ polynomial, then it has coordinate degree at most $D$. See [Kun24, § 1.3] or [O'D14, § 8.3] for further discussion.

## 1.4 Stability of Low Degree Algorithms

The key property of LDP and LCD algorithms for us is $L^2$ stability under input perturbations.

**Proposition 1.6** (Low Degree Stability). *Suppose $\mathcal{A}: \mathbf{R}^N \to \mathbf{R}^N$ is a deterministic algorithm with polynomial degree (resp. coordinate degree) $D$ and norm $\mathbf{E}\|\mathcal{A}(g)\|^2 \leq CN$. Then, for standard Normal r.v.s $g$ and $g'$ which are $(1-\varepsilon)$-correlated (resp. $(1-\varepsilon)$-resampled),*

$$\mathbf{E}\|\mathcal{A}(g) - \mathcal{A}(g')\|^2 \leq 2CD\varepsilon N, \tag{1.4}$$

*and thus for any $\eta > 0$,*

$$\mathbf{P}\big(\|\mathcal{A}(g) - \mathcal{A}(g')\| \geq 2\sqrt{\eta N}\big) \leq \frac{CD\varepsilon}{2\eta}. \tag{1.5}$$

*Proof:* We show (1.4) when $\mathcal{A}$ has coordinate degree $D$ and $(g, g')$ are $(1-\varepsilon)$-resampled. See e.g. [HS25b, Prop. 1.7] for the case where $\mathcal{A}$ is polynomial. In both cases, Markov's inequality gives (1.5).

We follow [Gam+22, Lem. 3.4]. Assume without loss of generality that $\mathbf{E}\|\mathcal{A}(g)\|^2 = 1$. Writing $\mathcal{A} = (\mathcal{A}_1, ..., \mathcal{A}_N)$, observe that for $g' \sim \mathcal{L}_{1-\varepsilon}(g)$,

$$\mathbf{E}\|\mathcal{A}(g) - \mathcal{A}(g')\|^2 = \mathbf{E}\|\mathcal{A}(g)\|^2 + \mathbf{E}\|\mathcal{A}(g')\|^2 - 2\mathbf{E}\langle\mathcal{A}(g), \mathcal{A}(g')\rangle = 2 - 2\mathbf{E}\langle\mathcal{A}(g), \mathcal{A}(g')\rangle.$$

By [O'D14, Exer. 8.18], we know that for each $\mathcal{A}_i \in V_{\leq D}$ we have

$$(1-\varepsilon)^D \mathbf{E}|\mathcal{A}_i(g)|^2 \leq \mathbf{E}[\mathcal{A}_i(g)\mathcal{A}_i(g')] \leq \mathbf{E}|\mathcal{A}_i(g)|^2.$$

Summing this over $i$ gives

$$(1-\varepsilon)^D \leq \mathbf{E}\langle\mathcal{A}(g), \mathcal{A}(g')\rangle \leq 1.$$

Combining this with the above, and using $1 - (1-\varepsilon)^D \leq \varepsilon D$, yields (1.4). □

**Remark 1.7.** Proposition 1.6 also holds for randomized algorithms, with exactly the same proof.

Next we introduce a family of locally improving algorithms, which will be useful later in showing the failure of randomized rounding. Below, fix a distance $r = O(1)$. Given a $\mathbf{R}^N$-valued $\mathcal{A}$, we can obtain a $\Sigma_N$-valued algorithm by first rounding $\mathcal{A}(g)$ into the solid hypercube $[-1, 1]^N$ and then picking the best corner of $\Sigma_N$ within constant distance of this output. Such a modification requires calculating the energy of at most $N^r$ additional points on $\Sigma_N$, and thus preserves e.g. any polynomial runtime bound. Since $r = O(1)$, this will also preserve stability. We formalize this construction as follows. Let clip: $\mathbf{R}^N \to [-1, 1]^N$ be the function which rounds $x \in \mathbf{R}^N$ into the cube $[-1, 1]^N$:

$$\mathsf{clip}(x)_i := \begin{cases} -1 & x_i \leq -1, \\ x_i & -1 < x_i < 1, \\ 1 & x_i \geq 1. \end{cases}$$

Note that clip is 1-Lipschitz with respect to the Euclidean norm.

**Definition 1.8.** Let $r > 0$ and $\mathcal{A}$ be an algorithm. Define the $[-1, 1]^N$-valued algorithm $\hat{\mathcal{A}}_r$ by

$$\hat{\mathcal{A}}_r(g) := \begin{cases} \underset{x' \in B_\Sigma(\mathsf{clip}(\mathcal{A}(g)), r)}{\arg\min} |\langle g, x' \rangle| & \text{if} \quad B_\Sigma(\mathsf{clip}(\mathcal{A}(g)), r) \neq \emptyset, \\ \mathsf{clip}(\mathcal{A}(g)) & \text{else.} \end{cases} \tag{1.6}$$

The next simply lemma shows that if $\mathcal{A}$ is stable, then $\hat{\mathcal{A}}_r$ is also stable.

**Lemma 1.9.** *Suppose $\mathcal{A}: \mathbf{R}^N \to \mathbf{R}^N$ is a deterministic algorithm with polynomial degree (resp. coordinate degree) $D$ and norm $\mathbf{E}\|\mathcal{A}(g)\|^2 \leq CN$. Then, for $r = O(1)$ and standard Normal r.v.s $g$ and $g'$ which are $(1-\varepsilon)$-correlated (resp. $(1-\varepsilon)$-resampled), $\hat{\mathcal{A}}_r$ as in Definition 1.8 satisfies:*

$$\mathbf{E}\|\hat{\mathcal{A}}_r(g) - \hat{\mathcal{A}}_r(g')\|^2 \leq 4CD\varepsilon N + 8r^2. \tag{1.7}$$

*Thus for any $\eta > 0$.*

$$\mathbf{P}\big(\|\hat{\mathcal{A}}_r(g) - \hat{\mathcal{A}}_r(g')\| \geq 2\sqrt{\eta N}\big) \leq \frac{CD\varepsilon}{\eta} + \frac{2r^2}{\eta N}. \tag{1.8}$$

*Proof*: Observe that by the triangle inequality, $\|\hat{\mathcal{A}}_r(g) - \hat{\mathcal{A}}_r(g')\|$ is bounded by

$$\|\hat{\mathcal{A}}_r(g) - \mathsf{clip}(\mathcal{A}(g))\| + \|\mathsf{clip}(\mathcal{A}(g)) - \mathsf{clip}(\mathcal{A}(g'))\| + \|\mathsf{clip}(\mathcal{A}(g')) - \hat{\mathcal{A}}_r(g')\|$$
$$\leq 2r + \|\mathcal{A}(g) - \mathcal{A}(g')\|.$$

This follows as clip is 1-Lipschitz and the corner-picking step in (1.6) only moves $\hat{\mathcal{A}}_r(g)$ from $\mathsf{clip}(\mathcal{A}(r))$ by at most $r$. By Jensen's inequality, squaring this gives

$$\|\hat{\mathcal{A}}_r(g) - \hat{\mathcal{A}}_r(g')\|^2 \leq 2\big(4r^2 + \|\mathcal{A}(g) - \mathcal{A}(g')\|^2\big).$$

Combining this with Proposition 1.6 gives (1.7), and (1.8) follows from Markov's inequality. □

Of course, our construction of $\hat{\mathcal{A}}_r$ is certainly never polynomial and does not preserve coordinate degree in a controllable way. However, because the rounding does not drastically alter the stability analysis, we are still able to show that for any $\mathbf{R}^N$-valued low coordinate degree algorithm $\mathcal{A}$ and $r = O(1)$, strong low degree hardness holds for $\hat{\mathcal{A}}_r$. Establishing the failure of $\hat{\mathcal{A}}_r$ will be a useful intermediate step towards the full proof of Theorem 1.4. (The same argument proves hardness when $\mathcal{A}$ is a low degree polynomial algorithm; this is omitted for brevity.)

## 2 Hardness for Low Degree Algorithms

In this section, we prove Theorem 1.3 and Theorem 1.4 – that is, we exhibit strong low degree hardness for both low polynomial degree and low coordinate degree algorithms.

Our argument utilizes what can be thought of as a "conditional" OGP. Previously, most OGP proofs identify a global obstruction: with high probability, there do not exist any tuples of good solutions to a family of correlated instances which are e.g. roughly the same distance apart. Here, however, we show a local obstruction; we condition on being able to solve a single instance and show that after a small change to the instance, it is unlikely that any solutions will exist close to the first one. This is an instance of the "brittleness" that makes the NPP challenging to solve; even small changes in the instance break the landscape geometry, so that even if solutions exist, there is no way to know where they will end up. Related strategies have been used recently in [AG24], [Ala24], [HS25b]. The first of these studied the Ising perceptron (in our terminology, the VBP with $d = \Theta(N)$), and deduced hardness of sampling from the fact that a *typical* solution to a given instance will disappear in a slightly correlated instance. For our result, it is crucial that *every* solution to a given NPP instance is likely to disappear when the input is slightly perturbed. On the other hand, this property is weaker than requiring all solutions to *simultaneously* disappear under such perturbation, as in a more standard "global" OGP argument. We note that Gamarnik and Kızıldağ showed in [GK23, Thm. 2.5] that sublinear energies do not exhibit a certain multi-OGP, suggesting that a sharp analysis via global OGP arguments may be challenging to implement.

Let us give a more detailed outline of our strategy, in the case of low coordinate degree. Let $E$ be an energy level, and assume $\mathcal{A}$ is a $\Sigma_N$-valued algorithm with coordinate degree at most $D \leq \tilde{o}(E)$. We choose suitable parameters $\eta \approx \frac{E}{N \log N}$ and $\varepsilon \approx \frac{\log(N/D)}{N}$, and aim to show that

$$\mathbf{P}(\mathcal{A}(g) \in S(E; g)) \to 0$$

as $N \to \infty$. To do so, we consider a $(1-\varepsilon)$-resampled pair $g, g'$ of NPP instances and proceed according to the following steps.

(a) For $\varepsilon$ small, $g$ and $g'$ have correlation close to 1. By Proposition 1.6, this implies that the outputs of an LCD algorithm $\mathcal{A}$ will be within distance $2\sqrt{\eta N}$ of each other with high probability.

(b) Since $\varepsilon \gg \frac{1}{N}$, we will have $g \neq g'$ with high probability, and we assume below that this holds. (This is the central difference between the LCD and LDP cases; in the latter there is no issue in taking $\varepsilon$ much smaller, which turns out to drastically affect the resulting hardness bound.)

(c) For $\eta$ small and fixed $\mathcal{A}(g)$, Lemma 2.4 shows using the Markov inequality that conditional on $g$ and the event $g' \neq g$, the perturbed instance $g'$ typically has no solutions within distance $2\sqrt{\eta N}$ of $\mathcal{A}(g)$. This is the conditional landscape obstruction we described above.

(d) Put together, these steps imply that it is unlikely for $\mathcal{A}$ to find solutions to *both* $g$ and $g'$ such that the stability guarantee of Proposition 1.6 holds. Using a positive correlation property (see Lemma 2.5), we conclude that $\mathcal{A}(g) \notin S(E; g)$ with high probability.

The above argument handles low coordinate degree algorithms without the randomized rounding step. To handle the latter, we observe that solutions to a given NPP instance are isolated with high probability. This implies that randomized rounding either changes only $O(1)$ coordinates (which preserves stability), or else injects too much randomness to preserve any chance of finding a good solution. For convenience, we pursue this extension only in the LCD case (as the main message of the LDP case is that polynomial degree is a poor proxy for complexity in the NPP).

## 2.1 Preliminary Estimates

We begin with some general estimates that will be utilized repeatedly throughout the proof. First, we bound the probability of a Normal r.v. being exponentially close to zero. We denote $2^x$ by $\exp_2(x)$.

**Lemma 2.1.** *Let $E, \sigma^2 > 0$, and suppose that conditionally on $\mu$, we have $Z \sim \mathcal{N}(\mu, \sigma^2)$. Then*

$$\mathbf{P}(|Z| \leq 2^{-E} \mid \mu) \leq \exp_2\left(-E - \frac{1}{2}\log_2(\sigma^2) + O(1)\right), \quad \forall \mu \in \mathbf{R}. \tag{2.1}$$

*Proof*: Observe that conditional on $\mu$, the distribution of $Z$ is bounded as

$$\varphi_{Z|\mu}(z) = \frac{1}{\sqrt{2\pi\sigma^2}} e^{-\frac{(z-\mu)^2}{2\sigma^2}} \leq (2\pi\sigma^2)^{-1/2}.$$

Integrating over $|z| \leq 2^{-E}$ then gives (2.1), via

$$\mathbf{P}(|Z| \leq 2^{-E}) = \int_{|z| \leq 2^{-E}} (2\pi\sigma^2)^{-1/2}\, \mathrm{d}z \leq 2^{-E - \frac{1}{2}\log_2(2\pi\sigma^2) + 1}. \qquad \square$$

Note that (2.1) is a decreasing function of $\sigma^2$. Thus, if there exists $\gamma$ with $\sigma^2 \geq \gamma > 0$, then (2.1) is bounded by $\exp_2(-E - \log_2(\gamma)/2 + O(1))$.

Next, we recall the following bound on the partial sums of binomial coefficients; this will be used for a first moment computation in Section 2.2.

**Lemma 2.2** (Chernoff-Hoeffding). *Suppose that $K \leq N/2$, and let $h(x) = -x\log_2(x) - (1-x)\log_2(x)$ be the binary entropy function. Then, for $p := K/N$,*

$$\sum_{k \leq K} \binom{N}{k} \leq \exp_2(Nh(p)) \leq \exp_2\left(2Np\log_2\left(\frac{1}{p}\right)\right).$$

*Proof*: For a $\mathrm{Bin}(N, p)$ random variable $S$, we have:

$$1 \geq \mathbf{P}(S \leq K) = \sum_{k \leq K} \binom{N}{k} p^k (1-p)^{N-k} \geq \sum_{k \leq K} \binom{N}{k} p^K (1-p)^{N-K}.$$

The last inequality follows by multiplying each term by $(p/(1-p))^{K-k} \leq 1$. Rearranging gives

$$\sum_{k \leq K} \binom{N}{k} \leq p^{-K}(1-p)^{-(N-K)}$$

$$= \exp_2(-K \log_2(p) - (N-K) \log_2(1-p))$$

$$= \exp_2\left(N \cdot \left(-\frac{K}{N} \log_2(p) - \left(\frac{N-K}{N}\right) \log_2(1-p)\right)\right)$$

$$= \exp_2(N \cdot (-p \log_2(p) - (1-p) \log_2(1-p))) = \exp_2(Nh(p)).$$

The final equality then follows from the bound $h(p) \leq 2p \log_2(1/p)$ for $p \leq 1/2$. □

Finally, we show that for any energy $E \leq O(N)$, there exists a choice of "distance" $\eta$ such that the term in the previous lemma is controlled.

**Lemma 2.3.** *For all $E \leq N$, there exist universal constants $C, C' > 0$ such that*

$$\eta := \frac{E}{C'N \log_2(CN/E)} \tag{2.2}$$

*satisfies $\eta \in (0, 1/2)$ and*

$$2\eta \log_2\left(\frac{1}{\eta}\right) < \frac{E}{4N}. \tag{2.3}$$

*In addition*
  (a) *if $E = \Theta(N)$, this $\eta$ is $\Theta(1)$;*
  (b) *if $E \leq o(N)$, this $\eta$ satisfies $\eta \gtrsim E/N \log_2(CN/E)$*

*Proof*: It is not hard to show that if $0 < \eta, E/N \ll 1$ and $-\eta \log \eta \sim E/N$, then we have $\eta \sim E/(N \log N/E)$. Picking suitable constants (e.g. $C = 8$ and $C' = 16$ suffice) in (2.2), it is easy to see that $\eta \in (0, 1/2)$ and (2.3) holds for all $E \leq N$. The resulting asymptotics follow immediately. □

## 2.2 Conditional Landscape Obstruction

We turn now to establishing the central conditional landscape obstruction of our argument. The idea is that for an $x \in \Sigma_N$ depending on $g$ and a related instance $g'$, the likelihood of any $x' \in \Sigma_N$ solving $g'$ is much smaller than the number of points within a neighborhood of $x$. Thus, even small changes to the instance destroy any solutions.

**Lemma 2.4.** *Let $(g, g')$ be a pair of either $(1-\varepsilon)$-correlated or $(1-\varepsilon)$-resampled instances. Let $x \in \Sigma_N$ be conditionally independent of $g'$ given $g$. Then for any $\eta \in (0, 1/2)$,*

$$\mathbf{P}_{g' \sim_{1-\varepsilon} g}\left(\begin{array}{l}\exists\, x' \in S(E; g') \text{ s.t.} \\ \|x - x'\| \leq 2\sqrt{\eta N}\end{array}\right) \leq \exp_2\left(-E - \frac{1}{2}\log_2(\varepsilon) + 2\eta \log_2\left(\frac{1}{\eta}\right)N + O(\log N)\right), \tag{2.4}$$

*and*

$$\mathbf{P}_{g' \sim \mathcal{L}_{1-\varepsilon}(g)}\left(\begin{array}{l}\exists\, x' \in S(E; g') \text{ s.t.} \\ \|x - x'\| \leq 2\sqrt{\eta N}\end{array}\middle|\, g \neq g'\right) \leq \exp_2\left(-E + 2\eta \log_2\left(\frac{1}{\eta}\right)N + O(1)\right). \tag{2.5}$$

*Proof*: Throughout, abbreviate $\eta_N := 2\sqrt{\eta N}$. We first show (2.4), by bounding the probability that

$$|B(x, \eta_N) \cap S(E; g')| = \sum_{x' \in B_\Sigma(x, \eta_N)} I\{x' \in S(E; g')\}$$

is nonzero. By Markov's inequality, this is upper bounded by

$$\mathbf{E}\left[\sum_{x' \in B_\Sigma(x, \eta_N)} I\{x' \in S(E; g')\}\right] = \mathbf{E}\left[\sum_{x' \in B_\Sigma(x, \eta_N)} \mathbf{E}[I\{x' \in S(E; g')\} \mid g]\right]$$
$$= \mathbf{E}\left[\sum_{x' \in B_\Sigma(x, \eta_N)} \mathbf{P}(|\langle g', x'\rangle| \leq 2^{-E} \mid g)\right]. \tag{2.6}$$

Note in particular that the range of this sum is independent of the inner probability, as $g'$ and $x$ are conditionally independent given $g$. To bound the number of terms in (2.6), let $k$ be the number of coordinates which differ between $x$ and $x'$, so that $\|x - x'\|^2 = 4k$. Thus $\|x - x'\| \leq 2\sqrt{\eta N}$ if and only if $k \leq N\eta < N/2$, so by Lemma 2.2,

$$|B_\Sigma(x, \eta_N)| = \sum_{k \leq N\eta} \binom{N}{k} \leq \exp_2(2\eta \log_2(1/\eta)N). \tag{2.7}$$

To bound the inner probability under $g' \sim_{1-\varepsilon} g$, fix any $x' \in \Sigma_N$ and write

$$g' = pg + \sqrt{1-p^2}\tilde{g}$$

for $p := 1 - \varepsilon$ and $\tilde{g}$ an independent copy of $g$. We know $\langle \tilde{g}, x'\rangle \sim \mathcal{N}(0, N)$, so this gives

$$\langle g', x'\rangle \mid g \sim \mathcal{N}(p\langle g, x'\rangle, (1-p^2)N).$$

This is nondegenerate for $(1-p^2)N \geq \varepsilon N > 0$; by Lemma 2.1, we get

$$\mathbf{P}_{g' \sim_{1-\varepsilon} g}(|\langle g', x'\rangle| \leq 2^{-E} \mid g) \leq \exp_2\left(-E - \frac{1}{2}\log_2(\varepsilon) + O(\log N)\right).$$

Using (2.7) to upper bound the number of terms in (2.6), summing this bound gives (2.4).

Alternatively, for (2.5), we know that if $g = g'$, then the $B(x, \eta_N) \cap S(E; g')$ will be nonempty if $x$ is chosen to be a solution to $g$; we thus condition on $g \neq g'$ throughout (2.6). To bound the corresponding inner term, again fix any $x' \in \Sigma_N$ and let $\tilde{g}$ be an independent copy of $g$. Let $J \subseteq [N]$ be a random subset where each $i \in J$ independently with probability $1 - \varepsilon$, so $g'$ can be represented as

$$g' = g_J + \tilde{g}_{\bar{J}}.$$

Conditional on $(g, J)$, we know that $\langle \tilde{g}_{\bar{J}}, x'\rangle$ is $\mathcal{N}(0, N - |J|)$ and $\langle g_J, x'\rangle$ is deterministic, so that

$$\langle g', x' \rangle \mid (g, J) \sim \mathcal{N}(\langle g_J, x' \rangle, N - |J|).$$

As $\{g \neq g'\} = \{|J| < N\}$, we have $N - |J| \geq 1$ conditional on $g \neq g'$. Thus, Lemma 2.1 gives

$$\mathbf{P}_{g' \sim \mathcal{L}_{1-\varepsilon}(g)}(|\langle g', x' \rangle| \leq 2^{-E} \mid g, g \neq g') \leq \exp_2(-E + O(1)),$$

and we conclude (2.5) as in the previous case. □

The following lemma shows a positive correlation property that enables us to avoid union-bounding over $\mathcal{A}$ finding solutions to correlated or resampled pairs of instances. Below, the set $S$ plays the role of the event $\mathcal{A}(g, \omega) \in S(E; g)$.

**Lemma 2.5** (Adapted from [HS25b, Lem. 2.7]). *Let $(g, g')$ be a pair of either $p$-correlated or $p$-resampled instances. Then for any set $S \subseteq \mathbf{R}^N$ and $p > 0$, with $q := \mathbf{P}(g \in S)$,*

$$\mathbf{P}_{g' \sim_p g}(g \in S, g' \in S) \geq q^2 \quad \text{and} \quad \mathbf{P}_{g' \sim \mathcal{L}_p(g)}(g \in S, g' \in S) \geq q^2.$$

*Proof*: In the first case, let $\tilde{g}, g^{(0)}, g^{(1)}$ be three i.i.d. copies of $g$, and observe that $g' \sim_p g$ are jointly representable as

$$g = \sqrt{1-\varepsilon}\tilde{g} + \sqrt{\varepsilon}g^{(0)}, \qquad g' = \sqrt{1-\varepsilon}\tilde{g} + \sqrt{\varepsilon}g^{(1)}.$$

Since $g, g'$ are conditionally i.i.d. given $\tilde{g}$, we have by Jensen's inequality that

$$\mathbf{P}_{g' \sim_p g}(g \in S, g' \in S) = \mathbf{E}[\mathbf{P}(g \in S, g' \in S \mid \tilde{g})] = \mathbf{E}[\mathbf{P}(g \in S \mid \tilde{g})^2] \geq \mathbf{E}[\mathbf{P}(g \in S \mid \tilde{g})]^2 = q^2.$$

Likewise, when $g' \sim \mathcal{L}_p(g)$, let $J$ be a random subset of $[N]$ where each $i \in J$ independently with probability $p$, so that $(g, g')$ are jointly representable as

$$g = \tilde{g}_J + g_{\bar{J}}^{(0)}, \qquad g' = \tilde{g}_J + g_{\bar{J}}^{(1)}.$$

Thus $g$ and $g'$ are conditionally i.i.d. given $(\tilde{g}, J)$, and we conclude in the same way. □

**Remark 2.6.** Note that Lemma 2.5 also holds when there is an auxiliary random seed $\omega$ shared across the instances. In this case the success event is a set $S \subseteq \mathbf{R}^N \times \Omega_N$, and we write

$$q = \mathbf{P}((g, \omega) \in S), \qquad Q = \mathbf{P}((g, \omega) \in S, (g', \omega) \in S),$$
$$q(\omega) = \mathbf{P}((g, \omega) \in S \mid \omega), \qquad Q(\omega) = \mathbf{P}((g, \omega) \in S, (g', \omega) \in S \mid \omega).$$

Lemma 2.5 shows that for any $\omega \in \Omega_N$, $Q(\omega) \geq q(\omega)^2$. Then, by Jensen's inequality,

$$Q = \mathbf{E}[q(\omega)] \geq \mathbf{E}[q(\omega)^2] \geq \mathbf{E}[q(\omega)]^2 = p^2.$$

Thus, in combination with Remark 1.7, the proofs of Theorem 1.3 and Theorem 1.4 hold without modification when $\mathcal{A}$ depends on an independent random seed $\omega$.

## 2.3 Hardness for LDP Algorithms

We first prove Theorem 1.3. Let $\mathcal{A}$ be a $\Sigma_N$-valued algorithm with satisfying (1.5) degree $D$; by Remark 2.6, assume without loss of generality that $\mathcal{A} = \mathcal{A}(g)$ is deterministic.

Consider a pair of $(1-\varepsilon)$-correlated instances $(g, g')$. Let $x := \mathcal{A}(g)$ and define the events

$$
\begin{aligned}
S_{\text{solve}} &:= \{\mathcal{A}(g) \in S(E; g), \mathcal{A}(g') \in S(E; g')\}, \\
S_{\text{stable}} &:= \{\|\mathcal{A}(g) - \mathcal{A}(g')\| \leq 2\sqrt{\eta N}\}, \\
S_{\text{cond}}(x) &:= \begin{Bmatrix} \nexists\ x' \in S(E; g') \text{ such that} \\ \|x - x'\| \leq 2\sqrt{\eta N} \end{Bmatrix}.
\end{aligned}
\tag{2.8}
$$

Intuitively, the first two events ask that the algorithm solves both instances and is stable, respectively. The last event, which depends on $x$, corresponds to the conditional landscape obstruction: for an $x$ depending only on $g$, there is no solution to $g'$ which is close to $x$.

**Lemma 2.7.** *For $x := \mathcal{A}(g)$, we have $S_{\text{solve}} \cap S_{\text{stable}} \cap S_{\text{cond}}(x) = \emptyset$.*

*Proof*: Suppose that $S_{\text{solve}}$ and $S_{\text{stable}}$ both occur. Letting $x := \mathcal{A}(g)$ (which only depends on $g$) and $x' := \mathcal{A}(g')$, we know $x' \in S(E; g')$ and is within distance $2\sqrt{\eta N}$ of $x$, contradicting $S_{\text{cond}}(x)$. □

*Proof of Theorem 1.3*: Let $p_{\text{solve}} := \mathbf{P}(\mathcal{A}(g) \in S(E; g))$ be the probability that $\mathcal{A}$ solves one instance. By Lemma 2.5, $\mathbf{P}_{g' \sim_{1-\varepsilon} g}(S_{\text{solve}}) \geq p_{\text{solve}}^2$.

In addition, let

$$p_{\text{cond}} := \max_{x \in \Sigma_N} 1 - \mathbf{P}_{g' \sim_{1-\varepsilon} g}(S_{\text{cond}}(x)). \qquad p_{\text{unstable}} := 1 - \mathbf{P}_{g' \sim_{1-\varepsilon} g}(S_{\text{stable}}),$$

Set $\varepsilon := 2^{-E/2}$ and set $\eta$ as in Lemma 2.3. By Lemma 2.4, for $N$ sufficiently large,

$$
\begin{aligned}
p_{\text{cond}} &\leq \exp_2\left(-E - \frac{1}{2}\log_2(\varepsilon) + 2\eta \log_2\left(\frac{1}{\eta}\right) N + O(\log N)\right) \\
&\leq \exp_2\left(-E + \frac{E}{4} + \frac{E}{4} + O(\log N)\right) = \exp_2\left(-\frac{E}{2} + O(\log N)\right) = o(1).
\end{aligned}
$$

Next, for $\log N \ll E \leq N$ and $D = o(\exp_2(E/4))$, (1.5) gives

$$
\begin{aligned}
p_{\text{unstable}} &\lesssim \frac{D\varepsilon}{\eta} \lesssim \frac{D \exp_2(-E/2) N \log_2(CN/E)}{E} \\
&\leq D \exp_2\left(-\frac{E}{2} + O(\log N)\right) \\
&\leq o(1) \cdot \exp_2\left(-\frac{E}{4} + O(\log N)\right) = o(1).
\end{aligned}
$$

That is, both $p_{\text{cond}}$ and $p_{\text{unstable}}$ vanish for large $N$.

To conclude, Lemma 2.7 gives (for $\mathbf{P} = \mathbf{P}_{g' \sim_{1-\varepsilon} g}$)

$$\mathbf{P}(S_{\text{solve}}) + \mathbf{P}(S_{\text{stable}}) + \mathbf{P}(S_{\text{cond}}(x)) \leq 2,$$

and rearranging yields

$$(p_{\text{solve}})^2 \leq p_{\text{unstable}} + (1 - \mathbf{P}(S_{\text{cond}}(x)) \leq p_{\text{unstable}} + p_{\text{cond}} = o(1). \qquad \square$$

## 2.4 Hardness for Non-Rounded LCD Algorithms

We now build towards the main result, Theorem 1.4. As a stepping stone, we first show strong low degree hardness for $\mathcal{A}$ with coordinate degree $D$ which are $\Sigma_N$-valued. We then extend to $O(1)$-close and rounded algorithms. Recalling Remark 2.6, we assume $\mathcal{A}$ is deterministic for convenience.

Consider a pair of $(1-\varepsilon)$-resampled instances $(g, g')$. Let $x := \mathcal{A}(g)$ and keep the definitions of $S_{\text{solve}}, S_{\text{stable}}, S_{\text{cond}}$ from (2.8). In addition, define

$$S_{\text{diff}} := \{g \neq g'\}.$$

**Lemma 2.8.** *For $x := \mathcal{A}(g)$, we have $S_{\text{diff}} \cap S_{\text{solve}} \cap S_{\text{stable}} \cap S_{\text{cond}}(x) = \emptyset$.*

*Proof*: This follows from Lemma 2.7, noting that the proof did not use that $g \neq g'$ almost surely. $\square$

As before, our proof follows from showing that for appropriate choices of $\varepsilon$ and $\eta$ (depending on $D$, $E$, and $N$), the events $S_{\text{stable}}$ and $S_{\text{cond}}(x)$ hold with high probability. However, when $S_{\text{diff}}$ fails, $S_{\text{cond}}(x)$ always fails. Thus, to ensure the appropriate probabilities vanish, we are required to choose $\varepsilon \gg 1/N$, which by (1.3) ensures $g \neq g'$ with high probability. Contrast this with Section 2.3, where $\varepsilon$ could be exponentially small in $E$. This restriction on $\varepsilon$ prevents us from showing hardness for algorithms with degree larger than the best possible level, as discussed in Section 1.2.

*Proof of Theorem 1.4, for $\Sigma_N$-valued algorithms **without randomized rounding***: Again let $p_{\text{solve}} := \mathbf{P}(\mathcal{A}(g) \in S(E; g))$ be the probability that $\mathcal{A}$ solves one instance. By Lemma 2.5, $\mathbf{P}_{g' \sim \mathcal{L}_{1-\varepsilon}(g)}(S_{\text{solve}}) \geq p_{\text{solve}}^2$. We now redefine $p_{\text{cond}}$ and $p_{\text{unstable}}$ via

$$p_{\text{cond}} := \max_{x \in \Sigma_N} 1 - \mathbf{P}_{g' \sim \mathcal{L}_{1-\varepsilon}(g)}(S_{\text{cond}}(x) \mid S_{\text{diff}}), \quad p_{\text{unstable}} := 1 - \mathbf{P}_{g' \sim \mathcal{L}_{1-\varepsilon}(g)}(S_{\text{stable}} \mid S_{\text{diff}}). \quad (2.9)$$

Setting $\eta$ as in Lemma 2.3, we have by Lemma 2.4 that

$$p_{\text{cond}} \leq \exp_2\left(-E + 2\eta \log_2\left(\frac{1}{\eta}\right) N + O(1)\right) \leq \exp_2\left(-\frac{3E}{4} + O(1)\right) = o(1).$$

Next, set $\varepsilon := \frac{\log_2(N/D)}{N}$. This clearly has $N\varepsilon \gg 1$, so

$$\mathbf{P}(S_{\text{diff}}) = 1 - (1-\varepsilon)^N \geq 1 - e^{-\varepsilon N} \to 1. \qquad (2.10)$$

By Proposition 1.6, we have for (a) $E = \delta N$ and $D = o(N)$ that

$$p_{\text{unstable}} \lesssim D\varepsilon = \frac{D}{N} \log_2\left(\frac{N}{D}\right) = o(1).$$

Likewise, for (b) $\log^2 N \ll E \ll N$ and $D = o(E/\log^2(CN/E))$, we get

$$p_{\text{unstable}} \lesssim \frac{D \log_2(N/D) \log_2(CN/E)}{E}$$

As $D$ is an upper bound on the maximum possible algorithm degree, we may increase $D$ without loss of generality in the analysis so that $D$ grows only slightly slower than $E$. Thus we assume henceforth that $D \geq E/\log_2^3(N/E)$, so that $N/D \leq N \log_2^3(N/E)/E$. This lets us bound

$$\log_2(N/D) \leq \log_2(N/E) + 3 \log_2 \log_2(N/E) \lesssim \log_2(CN/E),$$

which gives

$$p_{\text{unstable}} \lesssim \frac{D \log_2^2(CN/E)}{E} = o(1).$$

As before, these choices of $\varepsilon$ and $\eta$ ensure both $p_{\text{cond}}$ and $p_{\text{unstable}}$ vanish for large $N$ and arbitrary energy $\log^2 N \ll E \leq N$.

To conclude, for $x = \mathcal{A}(g)$, Lemma 2.8 implies $\mathbf{P}_{g' \sim \mathcal{L}_{1-\varepsilon}(g)}(S_{\text{solve}}, S_{\text{stable}}, S_{\text{cond}}(x) \mid S_{\text{diff}}) = 0$, so

$$\mathbf{P}(S_{\text{solve}}|S_{\text{diff}}) + \mathbf{P}(S_{\text{stable}}|S_{\text{diff}}) + \mathbf{P}(S_{\text{cond}}(x)|S_{\text{diff}}) \leq 2.$$

Thus, rearranging and multiplying by $\mathbf{P}(S_{\text{diff}})$ gives

$$\mathbf{P}(S_{\text{solve}} \text{ and } S_{\text{diff}}) \leq \mathbf{P}(S_{\text{diff}}) \cdot (p_{\text{unstable}} + p_{\text{cond}}) \leq p_{\text{unstable}} + p_{\text{cond}}.$$

Finally, adding $\mathbf{P}(S_{\text{solve}}, \neg S_{\text{diff}}) \leq 1 - P(S_{\text{diff}})$, which is $o(1)$ by our choice of $\varepsilon$, to both sides (so as to apply Lemma 2.5) lets us conclude

$$p_{\text{solve}}^2 \leq \mathbf{P}(S_{\text{solve}}) \leq p_{\text{unstable}} + p_{\text{cond}} + (1 - \mathbf{P}(S_{\text{diff}})) = o(1). \qquad \square$$

## 2.5 Locally Improving Algorithms

So far, we have established strong low degree hardness for both low degree polynomial and low coordinate degree algorithms which take values in $\Sigma_N$. Next we show that the last condition is not in fact as restrictive as it might appear, focusing on the LCD case. In the next step towards Theorem 1.4, we show that our preceding analysis extends to $O(1)$-distance perturbations of an LCD algorithm, thanks to the preservation of stability.

Throughout the below, fix a distance $r = O(1)$. We consider the event that the $\mathbf{R}^N$-valued algorithm $\mathcal{A}$ outputs a point close to a solution for an instance $g$:

$$S_{\text{close}}(r) = \begin{cases} \exists\, \hat{x} \in S(E;g) \text{ s.t.} \\ \text{clip}(\mathcal{A}(g)) \in B(\hat{x}, r) \end{cases} = \{B_{\Sigma}(\text{clip}(\mathcal{A}(g)), r) \neq \emptyset\}.$$

Equivalently, this means that the rounded algorithm $\hat{\mathcal{A}}_r$ defined in (1.6) solves the NPP.

**Proposition 2.9** (Hardness for Locally Improved LCD Algorithms)**.** *Let $g \sim \mathcal{N}(0, I_N)$ be a standard Normal random vector. Let $r > 0$ be an $N$-independent constant and $\mathcal{A}$ be any coordinate degree $D$ algorithm with $\mathbf{E}\|\mathcal{A}(g)\|^2 \lesssim N$. Assume that*
  *(a) if $E = \delta N$ for $\delta \in (0,1)$, then $D \leq o(N)$;*

(b) if $\omega(\log^2 N) \leq E \leq o(N)$, then $D \leq o(E/\log^2(CN/E))$.

Then $\mathbf{P}(S_{\text{close}}(r)) = \mathbf{P}(\hat{\mathcal{A}}_r(g) \in S(E;g)) = o(1)$.

*Proof of Proposition 2.9*: We maintain the setup of the proof of Theorem 1.4. Let $x := \hat{\mathcal{A}}_r(g)$ and define the events $S_{\text{diff}}, S_{\text{solve}}, S_{\text{stable}}$, and $S_{\text{cond}}(x)$ as in Section 2.4, replacing $\mathcal{A}$ with $\hat{\mathcal{A}}_r$. Let $p_{\text{solve}} := \mathbf{P}(\hat{\mathcal{A}}_r(g) \in S(E;g))$; letting $S$ be the event $\{\hat{\mathcal{A}}_r(g) \in S(E;g)\} = \{\text{clip}(\mathcal{A}(g)) \in S(E;g) + B(0,r)\}$, Lemma 2.5 ensures $\mathbf{P}_{g' \sim \mathcal{L}_{1-\varepsilon}(g)}(S_{\text{solve}}) \geq p_{\text{solve}}^2$. We keep the same definitions of $p_{\text{unstable}}$ and $p_{\text{cond}}$ as in (2.9), again replacing $\mathcal{A}$ with $\hat{\mathcal{A}}_r$. Thus, we choose $\varepsilon := \log_2(N/D)/N$, so that $\mathbf{P}(S_{\text{diff}}) \to 1$. Additionally, choose $\eta$ as Lemma 2.3, so that $p_{\text{cond}} = o(1)$.

It then suffices to show that $p_{\text{unstable}} = o(1)$. To see this, recall that when $E = \delta N$, our choice of $\eta$ is $\Theta(1)$, so $\frac{r^2}{\eta N} = o(1)$. In the sublinear case $\omega(\log^2 N) \leq E \leq o(N)$, we instead get

$$\eta N \gtrsim \frac{E}{N \log_2(CN/E)} \cdot N \geq \frac{E}{\log_2(CN/E)} = \omega(1),$$

since $E \gg \log^2 N$. Applying the properly modified Lemma 1.9 – knowing that the first term bounding $p_{\text{unstable}}$ is $o(1)$ with these choices of $\varepsilon$ and $\eta$ – we see that $\hat{p}_{\text{unstable}} = o(1)$, as desired. □

## 2.6 Truly Random Rounding

To complete the proof of Theorem 1.4, we need to show that randomized rounding cannot help solve the NPP. For this, we show in Lemma 2.12 below that if one has a subcube of $\Sigma_N$ with dimension growing slowly with $N$, then at most one of those points will be a solution. Then we show that randomized rounding fails with high probability whenever it changes a diverging number of coordinates. Together with the previous subsection which addresses $O(1)$ coordinate changes, this completes the proof. Below let $\mathcal{A}$ be an $\mathbf{R}^N$-valued algorithm and $\omega(\log^2 N) \leq E \leq N$.

We then define the deterministically rounded algorithm $\mathcal{A}^*(g) := \text{sign}(\mathcal{A}(g))$, which is the most likely single outcome of randomized rounding on $\mathcal{A}(g)$. Define $p_1(x), ..., p_N(x) \in [0, \tfrac{1}{2}]$ by

$$p_i(x) := \mathbf{P}(\text{round}(x)_i \neq \text{sign}(x)_i) = \frac{|x_i - \text{sign}(x_i)|}{2}. \tag{2.11}$$

(Here we suppress the dependence on $\vec{U} = (U_1, ..., U_N)$ from just above Theorem 1.4.)

**Lemma 2.10.** *Fix $x \in \mathbf{R}^N$. Draw $N$ coin flips $I_{x,i} \sim \text{Bern}(2p_i(x))$ as well as $N$ signs $S_i \sim \text{Unif}\{\pm 1\}$, all mutually independent; define the random variable $\tilde{x} \in \Sigma_N$ by*

$$\tilde{x}_i := S_i I_{x,i} + (1 - I_{x,i})\text{sign}(x)_i.$$

*Then $\tilde{x} \sim \text{round}(x)$.*

*Proof*: Let $x^* := \text{sign}(x)$. Conditioning on $I_{x,i}$, we can check that

$$\mathbf{P}(\tilde{x}_i \neq x_i^*) = 2p_i(x) \cdot \mathbf{P}(\tilde{x}_i = x_i \mid I_{x,i} = 1) + (1 - 2p_i(x)) \cdot \mathbf{P}(\tilde{x}_i \neq x_i^* \mid I_{x,i} = 0) = p_i(x).$$

Thus, $\mathbf{P}(\tilde{x}_i = x_i^*) = \mathbf{P}(\text{round}(x)_i = x_i^*)$. □

By Lemma 2.10, we can redefine round($x$) to be $\tilde{x}$ as constructed above without loss of generality. It thus makes sense to define $\tilde{\mathcal{A}}(g) := \text{round}(\mathcal{A}(g))$, which is now always $\Sigma_N$-valued. We have seen in Proposition 2.9 that when $\|\tilde{\mathcal{A}} - \mathcal{A}^*\| \leq O(1)$, low degree hardness will still apply. On the other hand, when $\|\tilde{\mathcal{A}} - \mathcal{A}^*\|_1 \geq \omega(1)$, any rounding scheme will introduce so much randomness that $\tilde{x}$ will effectively be a random point, which has a vanishing probability of being a solution due to the local sparsity of the solution set. (These two cases suffice as $\|\tilde{\mathcal{A}} - \mathcal{A}^*\| \leq \|\tilde{\mathcal{A}} - \mathcal{A}^*\|_1$.)

To see this, we first show that any randomized rounding scheme as in Lemma 2.10 which differs with high probability from $\mathcal{A}^*$ must resample a diverging number of coordinates.

**Lemma 2.11.** *Fix $x \in \mathbf{R}^N$, and let $p_1(x), ..., p_N(x)$ be defined as in (2.11). Then $\tilde{x} \neq x^*$ holds with probability $1 - o(1)$ if and only if $\sum_i p_i(x) = \omega(1)$. Moreover if $\sum_i p_i(x) = \omega(1)$, the number of coordinates $\tilde{x}$ is resampled in diverges in probability (i.e. exceeds any fixed $m < \infty$ with probability $1 - o(1)$).*

*Proof*: Recall that for $x \in [0, 1/2]$, we have $\log_2(1-x) = \Theta(x)$. Thus, as each coordinate of $x$ is rounded independently, we can compute

$$\mathbf{P}(\tilde{x} = x^*) = \prod_i (1 - p_i(x)) = \exp_2\left(\sum_i \log_2(1 - p_i(x))\right) \leq \exp_2\left(-\Theta\left(\sum_i p_i(x)\right)\right).$$

Thus, $\mathbf{P}(\tilde{x} = x^*) = o(1)$ if and only if $\sum_i p_i(x) = \omega(1)$.

Next, following the construction of $\tilde{x}$ in Lemma 2.10, let $E_i = \{I_{x,i} = 1\}$ be the event that $\tilde{x}_i$ is resampled from $\text{Unif}\{\pm 1\}$, independently of $x_i^*$. The $E_i$ are independent Bernoulli variables, so $\sum_i E_i$ has variance at most its mean. Hence the second claim follows by Chebyshev's inequality. $\square$

The following Lemma 2.12 takes advantage of this by showing that NPP solutions are isolated.

**Lemma 2.12** (Solutions Repel)**.** *Consider any distances $k = \Omega(1)$ and energy levels $E \gg k \log N$. With high probability, there are no pairs of distinct solutions $x, x' \in S(E; g)$ to an instance $g$ with $\|x - x'\| \leq 2\sqrt{k}$ (i.e., within $k$ sign flips of each other):*

$$\mathbf{P}\begin{pmatrix} \exists \, (x, x') \in S(E; g) \text{ s.t.} \\ \|x - x'\| \leq 2\sqrt{k}. \end{pmatrix} \leq \exp_2(-E + O(k \log N)) = o(1). \tag{2.12}$$

*Proof*: Consider any $x \neq x'$, and let $J \subseteq [N]$ denote the coordinates in which $x, x'$ differ. Then

$$x = x_{\overline{J}} + x_J, \qquad\qquad x' = x_{\overline{J}} - x_J.$$

Assuming both $x, x' \in S(E; g)$, we can expand the inequalities $-2^{-E} \leq \langle g, x \rangle, \langle g, x' \rangle \leq 2^{-E}$ into

$$-2^{-E} \leq \langle g, x_{\overline{J}} \rangle + \langle g, x_J \rangle \leq 2^{-E},$$
$$-2^{-E} \leq \langle g, x_{\overline{J}} \rangle - \langle g, x_J \rangle \leq 2^{-E}.$$

Multiplying the lower equation by $-1$ and adding the resulting inequalities gives $|\langle g, x_J \rangle| \leq 2^{-E}$.

Thus, finding pairs of distinct solutions within distance $2\sqrt{k}$ implies finding a subset $J \subseteq [N]$ of at most $k$ coordinates and $|J|$ signs $x_J$ such that $|\langle g_J, x_J \rangle| \leq 2^{-E}$. By [Ver18, Exer. 0.0.5], there are

$$\sum_{1 \leq k' \leq k} \binom{N}{k'} \leq \left(\frac{eN}{k}\right)^k \leq (eN)^k = 2^{O(k \log N)}$$

choices of such subsets, and at most $2^k$ choices of signs. Now, $\langle g_J, x_J \rangle \sim \mathcal{N}(0, |J|)$, and as $|J| \geq 1$, Lemma 2.1 and the following remark implies $\mathbf{P}(|\langle g_J, x_J \rangle| \leq 2^{-E}) \leq \exp_2(-E + O(1))$. Union bounding this over the $2^{O(k \log N)}$ possibilities gives (2.12). □

Next we deduce strong hardness for "truly randomized" algorithms which resample a diverging number of coordinates. Roughly, this holds because if enough coordinates are resampled, the resulting point is conditionally random within a subcube of dimension growing slowly with $\Sigma_N$, which by Lemma 2.12 can only contain a single solution.

**Theorem 2.13.** *Let $x = \mathcal{A}(g)$, and define $x^*, \tilde{x}$, etc., as previously. Moreover, assume that for any $x$ in the possible outputs of $\mathcal{A}$, we have $\sum_i p_i(x) = \omega(1)$. Then, for any $E \geq \omega(\log^2 N)$, we have*

$$\mathbf{P}(\tilde{\mathcal{A}}(g) \in S(E; g)) = \mathbf{P}(\tilde{x} \in S(E; g)) \leq o(1).$$

*Proof*: Following the characterization of $\tilde{x}$ in Lemma 2.10, let $K := \max(\log_2 N, \sum_i I_{x,i})$. By the assumptions on $\sum_i p_i(x)$ and Lemma 2.11, we know $K$, which is at least the number of coordinates which are resampled, is bounded as $1 \ll K \leq \log_2 N$, for any possible $x = \mathcal{A}(g)$. Now, let $J \subseteq [N]$ denote the set of the first $K$ coordinates to be resampled, so that $K = |J|$, and consider

$$\mathbf{P}(\tilde{x} \in S(E; g) \mid \tilde{x}_{\overline{J}}),$$

where we fix the coordinates outside of $J$ and let $\tilde{x}$ be uniformly sampled from a $K$-dimensional subcube of $\Sigma_N$. All such $\tilde{x}$ are within distance $2\sqrt{K}$ of each other, so by Lemma 2.12, the probability that there is more than one such $\tilde{x} \in S(E; g)$ is bounded by

$$\exp_2(-E + O(K \log N)) \leq \exp_2(-E + O(\log^2 N)) = o(1),$$

by assumption on $E$. Thus, the probability that any of the $\tilde{x}$ is in $S(E; g)$ is bounded by $2^{-K}$, whence

$$\mathbf{P}(\tilde{x} \in S(E; g)) = \mathbf{E}[\mathbf{P}(\tilde{x} \in S(E; g) \mid \tilde{x}_{\overline{J}})] \leq 2^{-K} \leq o(1). \qquad \square$$

Finally we combine the preceding results to deduce Theorem 1.4.

*Proof of Theorem 1.4*: Fix $r > 0$ an arbitrarily large integer (but fixed as $N$ grows). First, we know from Proposition 2.9 that $\mathbf{P}(\hat{\mathcal{A}}_r(g) \in S(E; g)) = o_{N \to \infty}(1)$ for any fixed $r$. Therefore by taking $r$ large slowly with $N$, it suffices to show that

$$\mathbf{P}(\tilde{\mathcal{A}}(g) \in S(E; g), \hat{\mathcal{A}}_r(g) \notin S(E; g)) \leq o_{r \to \infty}(1).$$

Indeed since $\|x - x^*\|_2 \leq \|x - x^*\|_1$, the left-hand probability tends to 0 uniformly as $r$ grows by Theorem 2.13. This completes the proof. □

**Acknowledgement.** Thanks to Brice Huang and Subhabrata Sen for comments and discussions.